\magnification\magstep1
\baselineskip15pt
\newread\AUX\immediate\openin\AUX=\jobname.aux
\def\ref#1{\expandafter\edef\csname#1\endcsname}
\ifeof\AUX\immediate\write16{\jobname.aux gibt es nicht!}\else
\input \jobname.aux
\fi\immediate\closein\AUX
\def\today{\number\day.~\ifcase\month\or
  Januar\or Februar\or M{\"a}rz\or April\or Mai\or Juni\or
  Juli\or August\or September\or Oktober\or November\or Dezember\fi
  \space\number\year}
\font\sevenex=cmex7
\scriptfont3=\sevenex
\font\fiveex=cmex10 scaled 500
\scriptscriptfont3=\fiveex

\def\phi{\varphi}
\def\epsilon{\varepsilon}
\def\theta{\vartheta}
\def\uauf{\lower1.7pt\hbox to 3pt{%
\vbox{\offinterlineskip
\hbox{\vbox to 8.5pt{\leaders\vrule width0.2pt\vfill}%
\kern-.3pt\hbox{\lams\char"76}\kern-0.3pt%
$\raise1pt\hbox{\lams\char"76}$}}\hfil}}
\def\cite#1{\expandafter\ifx\csname#1\endcsname\relax
{\bf?}\immediate\write16{#1 ist nicht definiert!}\else\csname#1\endcsname\fi}
\def\expandwrite#1#2{\edef\next{\write#1{#2}}\next}
\def\neverexpand{\noexpand\noexpand\noexpand}
\def\strip#1\ {}
\def\ncite#1{\expandafter\ifx\csname#1\endcsname\relax
{\bf?}\immediate\write16{#1 ist nicht definiert!}\else
\expandafter\expandafter\expandafter\strip\csname#1\endcsname\fi}
\newwrite\AUX
\immediate\openout\AUX=\jobname.aux
\newcount\Abschnitt\Abschnitt0
\def\beginsection#1. #2 \par{\advance\Abschnitt1%
\vskip0pt plus.10\vsize\penalty-250
\vskip0pt plus-.10\vsize\bigskip\vskip\parskip
\edef\TEST{\number\Abschnitt}
\expandafter\ifx\csname#1\endcsname\TEST\relax\else
\immediate\write16{#1 hat sich geaendert!}\fi
\expandwrite\AUX{\neverexpand\ref{#1}{\TEST}}
\leftline{\bf\number\Abschnitt. \ignorespaces#2}%
\nobreak\smallskip\noindent\SATZ1}
\def\Proof:{\par\noindent{\it Proof:}}
\def\Remark:{\ifdim\lastskip<\medskipamount\removelastskip\medskip\fi
\noindent{\bf Remark:}}
\def\Remarks:{\ifdim\lastskip<\medskipamount\removelastskip\medskip\fi
\noindent{\bf Remarks:}}
\def\Definition:{\ifdim\lastskip<\medskipamount\removelastskip\medskip\fi
\noindent{\bf Definition:}}
\def\Example:{\ifdim\lastskip<\medskipamount\removelastskip\medskip\fi
\noindent{\bf Example:}}
\newcount\SATZ\SATZ1
\def\proclaim #1. #2\par{\ifdim\lastskip<\medskipamount\removelastskip
\medskip\fi
\noindent{\bf#1.\ }{\it#2}\Par
\ifdim\lastskip<\medskipamount\removelastskip\goodbreak\medskip\fi}
\def\Aussage#1{%
\expandafter\def\csname#1\endcsname##1.{\ifx?##1?\relax\else
\edef\TEST{#1\penalty10000\ \number\Abschnitt.\number\SATZ}
\expandafter\ifx\csname##1\endcsname\TEST\relax\else
\immediate\write16{##1 hat sich geaendert!}\fi
\expandwrite\AUX{\neverexpand\ref{##1}{\TEST}}\fi
\proclaim {\number\Abschnitt.\number\SATZ. #1\global\advance\SATZ1}.}}
\Aussage{Theorem}
\Aussage{Proposition}
\Aussage{Corollary}
\Aussage{Lemma}
\font\la=lasy10
\def\strich{\hbox{$\vcenter{\hbox
to 1pt{\leaders\hrule height -0,2pt depth 0,6pt\hfil}}$}}
\def\dashedrightarrow{\hbox{%
\hbox to 0,5cm{\leaders\hbox to 2pt{\hfil\strich\hfil}\hfil}%
\kern-2pt\hbox{\la\char\string"29}}}

\def\Bindestrich{\penalty10000-\hskip0pt}
\let\_=\Bindestrich
\def\.{{\sfcode`.=1000.}}

\def\Rechts#1{\rlap{$\scriptstyle#1$}}
\def\Par{\par}
\def\:={\mathrel{\raise0,9pt\hbox{.}\kern-2,77779pt
\raise3pt\hbox{.}\kern-2,5pt=}}
\def\=:{\mathrel{=\kern-2,5pt\raise0,9pt\hbox{.}\kern-2,77779pt
\raise3pt\hbox{.}}} 

\def\pfeil{\rightarrow}

\def\pf#1{\buildrel#1\over\rightarrow}

\def\Ugleich{\hbox{$\cup$\kern.5pt\vrule depth -0.5pt}}
\def\|#1|{\mathop{\rm#1}\nolimits}
\def\<{\langle}
\def\>{\rangle}
\let\Times=\times
\def\times{\mathop{\Times}}
\let\Otimes=\otimes
\def\otimes{\mathop{\Otimes}}
\catcode`\@=11
\def\hex#1{\ifcase#1 0\or1\or2\or3\or4\or5\or6\or7\or8\or9\or A\or B\or
C\or D\or E\or F\else\message{Warnung: Setze hex#1=0}0\fi}
\def\fontdef#1:#2,#3,#4.{%
\alloc@8\fam\chardef\sixt@@n\FAM
\ifx!#2!\else\expandafter\font\csname text#1\endcsname=#2
\textfont\the\FAM=\csname text#1\endcsname\fi
\ifx!#3!\else\expandafter\font\csname script#1\endcsname=#3
\scriptfont\the\FAM=\csname script#1\endcsname\fi
\ifx!#4!\else\expandafter\font\csname scriptscript#1\endcsname=#4
\scriptscriptfont\the\FAM=\csname scriptscript#1\endcsname\fi
\expandafter\edef\csname #1\endcsname{\fam\the\FAM\csname text#1\endcsname}
\expandafter\edef\csname hex#1fam\endcsname{\hex\FAM}}
\catcode`\@=12 

\fontdef Ss:cmss10,,.
\fontdef Fr:eufm10,eufm7,eufm5.


\def\fQ{{\Fr Q}}

\fontdef bbb:msbm10,msbm7,msbm5.
\fontdef mbf:cmmib10,cmmib7,.
\fontdef msa:msam10,msam7,msam5.
\def\CC{{\bbb C}}

\def\NN{{\bbb N}}

\def\ZZ{{\bbb Z}}
\def\cD{{\cal D}}
\def\cE{{\cal E}}
\def\cI{{\cal I}}
\def\cP{{\cal P}}

\mathchardef\leer=\string"0\hexbbbfam3F
\mathchardef\subsetneq=\string"3\hexbbbfam24
\mathchardef\semidir=\string"2\hexbbbfam6E
\mathchardef\dirsemi=\string"2\hexbbbfam6F
\mathchardef\haken=\string"2\hexmsafam78
\mathchardef\auf=\string"3\hexmsafam10
\let\OL=\overline
\def\overline#1{{\hskip1pt\OL{\hskip-1pt#1\hskip-1pt}\hskip1pt}}

\def\fQ{{\overline{f}}}

\def\Pq{{\overline{P}}}

%
\abovedisplayskip 9.0pt plus 3.0pt minus 3.0pt
\belowdisplayskip 9.0pt plus 3.0pt minus 3.0pt
\newdimen\Grenze\Grenze2\parindent\advance\Grenze1em
\newdimen\Breite
\newbox\DpBox
\def\NewDisplay#1$${\Breite\hsize\advance\Breite-\hangindent
\setbox\DpBox=\hbox{\hskip2\parindent$\displaystyle{#1}$}%
\ifnum\predisplaysize<\Grenze\abovedisplayskip\abovedisplayshortskip
\belowdisplayskip\belowdisplayshortskip\fi
\global\futurelet\nexttok\WEITER}
\def\WEITER{\ifx\nexttok\qed\expandafter\leftQEDdisplay
\else\leftdisplay\fi}
\def\leftdisplay{\hskip-\hangindent\leftline{\box\DpBox}$$}
\def\leftQEDdisplay{\hskip-\hangindent
\line{\copy\DpBox\hfill\lower\dp\DpBox\copy\QEDbox}%
\belowdisplayskip0pt$$\bigskip\let\nexttok=}
\everydisplay{\NewDisplay}
\newbox\QEDbox
\newbox\nichts\setbox\nichts=\vbox{}\wd\nichts=2mm\ht\nichts=2mm
\setbox\QEDbox=\hbox{\vrule\vbox{\hrule\copy\nichts\hrule}\vrule}
\def\qed{\leavevmode\unskip\hfil\null\nobreak\hfill\copy\QEDbox\medbreak}
\newdimen\HIindent
\newbox\HIbox
\def\setHI#1{\setbox\HIbox=\hbox{#1}\HIindent=\wd\HIbox}
\def\HI#1{\par\hangindent\HIindent\hangafter=0\noindent\leavevmode
\llap{\hbox to\HIindent{#1\hfil}}\ignorespaces}

\newdimen\maxSpalbr
\newdimen\altSpalbr

\def\beginrefs{%
\expandafter\ifx\csname Spaltenbreite\endcsname\relax
\def\Spaltenbreite{1cm}\immediate\write16{Spaltenbreite undefiniert!}\fi
\expandafter\altSpalbr\Spaltenbreite
\maxSpalbr0pt
\def\L|Abk:##1|Sig:##2|Au:##3|Tit:##4|Zs:##5|Bd:##6|S:##7|J:##8||{%
\edef\TEST{[##2]}
\expandafter\ifx\csname##1\endcsname\TEST\relax\else
\immediate\write16{##1 hat sich geaendert!}\fi
\expandwrite\AUX{\neverexpand\ref{##1}{\TEST}}
\setHI{[##2]\ }
\ifnum\HIindent>\maxSpalbr\maxSpalbr\HIindent\fi
\ifnum\HIindent<\altSpalbr\HIindent\altSpalbr\fi
\HI{[##2]}
\ifx-##3\relax\else{##3}: \fi
\ifx-##4\relax\else{##4}{\sfcode`.=3000.} \fi
\ifx-##5\relax\else{\it ##5\/} \fi
\ifx-##6\relax\else{\bf ##6} \fi
\ifx-##8\relax\else({##8})\fi
\ifx-##7\relax\else, {##7}\fi\Par}
\def\B|Abk:##1|Sig:##2|Au:##3|Tit:##4|Reihe:##5|Verlag:##6|Ort:##7|J:##8||{%
\edef\TEST{[##2]}
\expandafter\ifx\csname##1\endcsname\TEST\relax\else
\immediate\write16{##1 hat sich geaendert!}\fi
\expandwrite\AUX{\neverexpand\ref{##1}{\TEST}}
\setHI{[##2]\ }
\ifnum\HIindent>\maxSpalbr\maxSpalbr\HIindent\fi
\ifnum\HIindent<\altSpalbr\HIindent\altSpalbr\fi
\HI{[##2]}
\ifx-##3\relax\else{##3}: \fi
\ifx-##4\relax\else{##4}{\sfcode`.=3000.} \fi
\ifx-##5\relax\else{(##5)} \fi
\ifx-##7\relax\else{##7:} \fi
\ifx-##6\relax\else{##6}\fi
\ifx-##8\relax\else{ ##8}\fi\Par}
\def\Pr|Abk:##1|Sig:##2|Au:##3|Artikel:##4|Titel:##5|Hgr:##6|Reihe:{%
\edef\TEST{[##2]}
\expandafter\ifx\csname##1\endcsname\TEST\relax\else
\immediate\write16{##1 hat sich geaendert!}\fi
\expandwrite\AUX{\neverexpand\ref{##1}{\TEST}}
\setHI{[##2]\ }
\ifnum\HIindent>\maxSpalbr\maxSpalbr\HIindent\fi
\ifnum\HIindent<\altSpalbr\HIindent\altSpalbr\fi
\HI{[##2]}
\ifx-##3\relax\else{##3}: \fi
\ifx-##4\relax\else{##4}{\sfcode`.=3000.} \fi
\ifx-##5\relax\else{In: \it ##5}. \fi
\ifx-##6\relax\else{(##6)} \fi\PrII}
\def\PrII##1|Bd:##2|Verlag:##3|Ort:##4|S:##5|J:##6||{%
\ifx-##1\relax\else{##1} \fi
\ifx-##2\relax\else{\bf ##2}, \fi
\ifx-##4\relax\else{##4:} \fi
\ifx-##3\relax\else{##3} \fi
\ifx-##6\relax\else{##6}\fi
\ifx-##5\relax\else{, ##5}\fi\Par}
\bgroup
\baselineskip12pt
\parskip2.5pt plus 1pt
\hyphenation{Hei-del-berg}
\sfcode`.=1000
\beginsection References. References

}
\def\endrefs{%
\expandwrite\AUX{\neverexpand\ref{Spaltenbreite}{\the\maxSpalbr}}
\ifnum\maxSpalbr=\altSpalbr\relax\else
\immediate\write16{Spaltenbreite hat sich geaendert!}\fi
\egroup}

\Aussage{Example}
\fontdef CMEX:cmex10 scaled 1200,,.
\mathchardef\Sum="1\hexCMEXfam50
\def\rho{\varrho}

\font\BF=cmbx10 scaled \magstep2
\font\CSC=cmcsc10 
{\baselineskip2\baselineskip\rightskip0pt plus 5truecm\hyphenpenalty10000
\leavevmode\vskip0truecm\noindent
\BF Difference Equations and Symmetric Polynomials Defined by Their Zeros
}
\vskip1truecm
\leftline{{\CSC Friedrich Knop \& Siddhartha Sahi}%
\footnote*{\rm The authors were partially supported by NSF grants.}}
\leftline{Department of Mathematics, Rutgers University, 
New Brunswick NJ 08903, USA}
\leftline{knop@math.rutgers.edu, sahi@math.rutgers.edu}
\vskip1truecm

\def \lam {\lambda}

\def \sig {\sigma}
\def \eps {\varepsilon}

\def \Plr {P_\lam^{\rho}}

\def \Pr {P_\lam^{(1/r)}}
\def \Pla {P_\lam^{(\alpha)}}
\def \Pma {P_\mu^{(\alpha)}}
\def \Plrd {P_\lam^{r\delta}}
\def \Pmrd {P_\mu^{r\delta}}

\def \Jla {J_\lam^{(\alpha)}}

\def \Jlrd {J_\lam^{r\delta}}

\def \D {{\cal D}(t;r)}
\def \E {{\cal E}}

\def\Sd-1nrd{S_{d-1}^n(r\delta)}

\def\Sd-1nr{S_{d-1}^n(\rho)}
\beginsection Introduction. Introduction

In this paper, we are starting a systematic analysis of a class of symmetric
polynomials which, in full generality, has been introduced in \cite{Sahi}.
The main features of these functions are that they are defined by vanishing
conditions and that they are non\_homogeneous. They depend on several
parameters but we are studying mainly a certain subfamily which is indexed
by one parameter $r$. As a special case, we obtain for $r=1$ the factorial
Schur functions discovered by Biedenharn and Louck \cite{BL}.

Our main result is that for general $r$ these functions are eigenvalues of
difference operators, which are difference analogues of the
Sekiguchi\_Debiard differential operators. Thus the functions under
investigation are non\_homogeneous variants of Jack polynomials.

More precisely, let $\Lambda$ be the set of partitions of length $n$,
i.e., sequences of integers $(\lambda_i)$ with
$\lambda_1\ge\ldots\ge\lambda_n\ge0$. The degree $|\lambda|$ of a
partition $\lambda$ is the sum of its parts. Choose a vector
$\rho\in\CC^n$ which has to satisfy a mild condition. Then for every
$\lambda\in\Lambda$ there is (up to a constant) a unique symmetric
polynomial $P_\lambda$ of degree at most $d$ which satisfies the
following vanishing condition:
$$
\hbox{$P_\lambda(\mu+\rho)=0$ for
all partitions $\mu$ with $|\mu|\le|\lambda|$ and $\mu\ne\lambda$.}
$$
This kind of vanishing comes up in the study of invariant
differential operators and Capelli type identities on multiplicity free
spaces and has been, in special cases, observed by other
authors (e.g. \cite{HU}, \cite{Ok}).  

In full generality, we have basically only one result (beyond their
existence) about the polynomials $P_\lambda$, namely two explicit formulas
for $P_\lambda$ when $\lambda=1^k$. From then on, we are only considering
$\rho=r\delta$, where $r\in\CC$ and $\delta=(n-1,n-2,\ldots,1,0)$.

We prove that these $P_\lambda$ are simultaneous eigenfunctions of $n$
commuting {\it difference\/} operators. On the highest homogeneous part of a
polynomial, these difference operators act like well known differential
operators: the Sekiguchi\_Debiard operators. The eigenfunctions of those
are the Jack polynomials. This has as immediate consequence that the top
homogeneous part of $P_\lambda$ is a Jack polynomial.

In the later sections, we draw several conclusions from the difference
equations. As an application to the ``classical'' theory we give a new
proof of the Pieri rule for Jack polynomials using the polynomials
$P_\lambda$.

We conclude with a brief discussion of the ``integral'' form $J_\lam$
which in the homogeneous case, is a rescaling of the $P_\lam$ by a 
certain hooklength factor.  It turns out that the corresponding 
{\it inhomogeneous\/} polynomial seems to have integrality and positivity 
properties which generalize a conjecture of Macdonald for the homogeneous 
case. In this connection, we have recently proved some integrality and 
positivity results which we shall report on elsewhere.

\medskip\noindent
{\bf Acknowledgment:} We would like to thank G.~Olshanski for sending us
his paper \cite{Olsh}. It initiated most of the research to the
present paper. Furthermore, we would like to thank A.~Zelevinski telling
us about Olshanski's work.

\beginsection construction. The basic construction

The results of this section are essentially in \cite{Sahi}, however
in order to keep the development self-contained we give a quick
rederivation.

Let us write $S(n,d)\subset\ZZ^n$ for the set of partitions
$\lambda_1\ge\ldots\lam_n\ge0$ with $|\lam|:=\sum\lambda_i=d$.
We say that $\rho\in\CC^n$ is {\it dominant\/}
if $\rho_i-\rho_j\ne-1,-2,-3,\ldots$ for all $i<j$. Slightly weakening this
condition, we define $\rho$ to be
{\it$d$\_dominant} if $\rho_i-\rho_j\ne-1,-2,-3,\ldots,-\left\lfloor{d\over
i}\right\rfloor$ for all $i<j$ where $d\in\NN$.

\Theorem TK1. For any $d\in\NN$ and $\rho\in\CC^n$ put
$M:=S(n,d)+\rho\subseteq\CC^n$. Assume, $\rho$ is $d$\_dominant. Then for
every map $\fQ:M\pfeil\CC$ there is a unique symmetric polynomial $f$
of degree at most $d$ such that $f|_M=\fQ$.

\Proof: For any partition $\lam\in\ZZ^n$
let $m_\lam$ be the corresponding monomial symmetric function in $n$
variables. If we express an arbitrary symmetric function of degree $\le d$ in
terms of $m_\lam$, then the interpolation problem gives a
{\it square\/} system of linear equations for the coefficients. Hence
existence implies uniqueness.

To show existence, we argue by induction on $n+d$. The
case $n=0$ is vacuous, so we assume $n\ge 1$.

To any $\lambda\in S(n-1,d)$ we can append a zero and obtain a
partition $\lambda,0\in S(n,d)$. This way, we can define map
$g=\sum a_\lam m_\lam\mapsto g^+=\sum a_\lam m_{\lam,0}$. It is an 
injective map from symmetric functions
in $n-1$ variables to symmetric functions in $n$ variables.
It has the property that $g^+$ has the same degree as $g$, and 
$g^+(x_1,\ldots,x_{n-1},0)=g(x_1,\ldots,x_{n-1})$.

We will construct $f$ as a function of the form 
$$
f(x)=g^+(x_1-\rho_n,\ldots, x_n-\rho_n)+
\left[\prod_{i=1}^n(x_i-\rho_n)\right]h(x_1-1,\ldots,x_n-1)
$$
First, let us consider the set $M_0$ of all points $x=\lam+\rho\in M$ with
$\lam_n=0$. Since $x_n-\rho_n=0$, the first term equals
$g(x_1-\rho_n, \ldots,x_{n-1}-\rho_n)$ and the second term
vanishes. If $x$ runs through $M_0$ then
$x'=(x_1-\rho_n,\ldots,x_{n-1}-\rho_n)$ runs through $S(n-1,d)+\rho'$,
where $\rho':= (\rho_1-\rho_n,\ldots, \rho_{n-1}-\rho_n)$ which is also
$d$\_dominant. By induction we can find $g$ of degree $\le d$ with
$f(x)=g(x')=\fQ(x)$ for all $x\in M_0$.

Next, we consider the points $x\in M\setminus M_0$, i.e.,
$x=\lam+\rho\in M$ with $\lam_n>0$. These exist only if $d\ge n$. As $x$
runs through these points $(x_1-1,\ldots,x_n-1)$ will run through
$S(n,d-n)+\rho$. Since $\lfloor d/i\rfloor\ge\lam_i\ge\lam_n>0$ and since
$\rho$ is $d$\_dominant, each of the factors $x_i-\rho_n=
\lam_i+\rho_i-\rho_n$ is non-zero. By induction, we can find $h$ of
degree $\le d-n$ such that $h$ has prescribed values at $M\setminus M_0$.
\qed

We assume from now on that $\rho$ is dominant.  With the theorem, we
are going to define interpolation polynomials. To get the most
convenient normalization, we have to introduce some more notation:
Recall that a partition $\lambda$ can be represented by its {\it
diagram}, i.e., the set of all lattice points (called boxes)
$(i,j)\in\ZZ^2$ with $1\le i\le n$ and $1\le j\le\lambda_i$. The dual
partition $\lambda'$ is the one with the transposed diagram. Now, for
every box $s$ we define the {\it $\rho$-hooklength} to be
$c_\lambda^\rho(s):=(\lambda_i-j+1)+(\rho_i-\rho_{\lambda_j'})$ and
$c_\lam^\rho:=\prod_{s\in\lambda}c_\lambda^\rho(s)$.

\Definition: For any partition $\lambda\in S(n,d)$ let $P_\lambda^\rho$
be the unique polynomial in $n$ variables such that
\item{(1)} $P_\lambda^\rho$ is symmetric;
\item{(2)} $\|deg|P_\lambda^\rho=d$;
\item{(3)} $P_\lambda^\rho(\mu+\rho)=0$ for all $\mu\in S(n,d)$,
$\mu\ne\lambda$;
\item{(4)} $P_\lambda^\rho(\lambda+\rho)=c_\lambda^\rho$.

\medskip\noindent
The normalization condition (4) is motivated by the following
theorem. In fact, we could replace (4) by it.

\Theorem T3. Let $\Plr=\sum_{\mu:|\mu|\le|\lam|} u_{\lam\mu}^\rho
m_\mu$ be the expression in terms of monomial symmetric functions. Then
$u_{\lam\lam}^\rho=1$.

\Proof: 
We proceed by induction on $n+|\lam|$. As in the proof of \cite{TK1}
we express
$$
\Plr=g^+(x_1-\rho_n,\ldots, x_n-\rho_n)+
\left[\prod_{i=1}^n(x_i-\rho_n)\right]h(x_1-1,\ldots,x_n-1)
$$
First assume $\lam_n=0$. Put $\nu:=(\lam_1,\ldots,\lam_{n-1})$
and $\rho':= (\rho_1-\rho_n,\ldots,\rho_{n-1}-\rho_n)$. Then
\cite{TK1} implies $g=a P^{\rho'}_\nu$ with $a\in\CC^*$. Now, we
compare values at $x=\lambda+\rho$. Since
$c_\lam^\rho=c^{\rho'}_\nu$ we obtain $a=1$ and the
assertion follows by induction.

Next, suppose $\lam_n>0$. Then \cite{TK1} implies $g=0$ and
$h=a P_\nu^\rho(x_1-1,\ldots,x_n-1)$ where
$\nu:=(\lam_1-1,\ldots,\lam_n-1)$ and $a\in\CC^*$. Again, we compare
values at $x=\lambda+\rho$.  The linear factors are just the
$\rho$-hooklengths for the first column of $\lam$. Thus, $a=1$ and the
assertion follows by induction.\qed

Additionally, we got the following reduction formula:

\Corollary. Assume $\lambda$ is a partition with $\lambda_n>0$ and let
$\lambda^*\:=(\lambda_1-1,\ldots,\lambda_n-1)$. Then
$P_\lambda^\rho=\prod_i(x_i-\rho_n)P_{\lambda^*}^\rho(x_1-1,\ldots,x_n-1)$.

\beginsection spcl. Special cases

We don't know an explicit formula for $\Plr$ in general but several
special cases are known.

For arbitrary $\rho$ we have only a formula for $\lam=1^k$. This is
the partition with $k$ ones and $(n-k)$ zeros.  The functions
$P^\rho_{1^k}$ are important since they are analogues of the
elementary symmetric functions. In particular, they generate the
symmetric polynomials as a ring. Actually, we have {\it two\/}
formulas for them.

Recall that the elementary symmetric function $e_j(x)$ and the 
complete symmetric function $h_j(y)$ are the coefficients of
$t^j$ in the expansions of $E(x,t)= \prod_i(1+tx_i)$ 
and $H(y,t)=\prod_i(1-ty_i)^{-1}$ respectively.

\Proposition P2. Let $\rho$ be dominant and $1\le k\le n$. Then
$$
P_{1^k}^\rho=
\sum_{j=0}^k(-1)^{k-j}h_{k-j}(\rho_k,\ldots,\rho_n)e_j(x)=
\sum_{i_1<\ldots<i_k}\prod_{j=1}^k(x_{i_j}-\rho_{i_j+k-j}).
$$

\Proof: Denote the first expression by $P'$, the second by $P''$. We
are going to show that they both satisfy the definition of
$P_{1^k}^\rho$. Both have certainly the right degree and $m_{1^k}$ has
the right coefficient.

For the vanishing condition (3), let $x=\mu+\rho$ with $|\mu|\le k$
and $\mu\ne 1^k$. This forces $\mu_k=\ldots=\mu_n=0$ and
$x_k=\rho_k,\ldots,x_n=\rho_n$.  Observe that $P'$ is precisely the
coefficient of $t^k$ in the power series expansion of
$\prod_{i=1}^n(1+tx_i)/ \prod_{i=k}^n(1+t\rho_i)$.  Evaluated at $x$,
this quotient becomes a {\it polynomial\/} of degree $<k$, and its
$k$-th coefficient $P'(x)$ vanishes. As for $P''$, the index $i_k$ in its
definition is at least $k$. Hence the factors for $j=k$ vanish
at $x$ which shows $P''(x)=0$.

Finally, we have to show symmetry. This is trivial for $P'$ but not
quite for $P''$. First let $n=2$. Then
$$
P''_{1^1}=(x_1-\rho_1)+(x_2-\rho_2);\quad
P''_{1^2}=(x_1-\rho_2)(x_2-\rho_2) $$
which are certainly symmetric. Now let $n\ge3$.
To make the dependence on $\rho$ and $k$ visible, we write
$P''=P''_k(x;\rho)$. Furthermore, let $x'$, $\rho'$ (resp. $x''$, $\rho''$)
equal $x$, $\rho$ where we dropped the last (resp. first) component.
If we break the defining sum for $P''$ up according
$i_k<n$ or $i_k=n$ we get
$$
P''_k(x;\rho)=P''_k(x';\rho')+(x_n-\rho_n)P''_{k-1}(x';\rho'').
$$
By induction we see that $P''$ is symmetric in
$x_1,\ldots,x_{n-1}$. If we break the sum up according $i_1=1$ or not
we obtain
$$
P''_k(x;\rho)=P''_k(x'';\rho'')+(x_1-\rho_k)P''_{k-1}(x'';\rho'').
$$
This shows that $P''$ is symmetric in $x_2,\ldots,x_n$ as well.\qed

\Remarks: For $\rho=r(n-1,\ldots,1,0)$, the expression $P'$ is
essentially due to Wallach while that for $P''$ can be traced back to
Capelli. The equality $P'=P''$ can be also proved directly by using
the polynomials $e_k(x/y)$ of \cite{Macd3} p.58.
\medskip

For the rest of the paper we specialize to $\rho$ of the form
$r\delta$ where $r$ is a complex number or just an indeterminate and
$\delta:=(n-1,\ldots,1,0)$. The dominance of $\rho$ means that
$r\ne -p/q$ where $p,q$ are integers such that $p,q\ge 1$, and $q<n$.
We shall assume this from now on.

\let\UL=\underline
\def\underline#1{{\hskip.5pt\UL{\hskip-.5pt#1\hskip-.5pt}\hskip.5pt}}
\def\Fall#1{^{\underline{#1}}}
\def\fall#1{_{\underline{#1}}}

First we treat the case $r=0$. For this we introduce the {\it falling
factorial polynomials\/} $x\Fall m:=x(x-1)\ldots(x-m+1)$. The
factorial monomial symmetric functions $m\fall\lam$ are
obtained by replacing each monomial 
$x_1^{l_1}x_2^{l_2} \ldots x_n^{l_n}$ in $m_\lam$ by the corresponding 
factorial monomial $x_1\Fall{l_1}x_2\Fall{l_2}\ldots x_n\Fall{l_n}$.
The following is obvious.

\Proposition P4. For $r=0$, we have $P_\lam^0= m\fall\lam$.

For $r=1$ we get the factorial Schur functions. (See \cite{BL}, \cite{Macd2},
and \cite{Olsh}.) To define them we write $a_\delta(x)$
for the Vandermonde determinant $\|det|\big(x_i^{\delta_j}\big)=
\prod_{i<j}(x_i-x_j)$. Then the next result seems to be due to Okounkov
\cite{Ok}.

\Proposition P5. For $r=1$, we have 
$$
P_\lambda^\delta(x)=
{1\over a_\delta(x)}\|det|\big(x_i\Fall{\lambda_j+\delta_j}\big).
$$

\Proof: 
Since $\|det|\big(x_i\Fall{\lambda_j+\delta_j}\big)$ is a skew
symmetric polynomial, its quotient by $a_\delta$ is a symmetric
polynomial which is easily seen to have degree $|\lam|$. Now let
$\mu\ne\lam$ and $|\mu|\le|\lam|$. Since $a_\delta(\mu+\delta)\ne0$
for any partition $\mu$, it remains only to prove the vanishing of
$\|det|\big[(\mu_i+\delta_i)\Fall{\lambda_j+\delta_j}\big]=\sum_\sig
(-1)^\sigma\prod_i(\mu_{\sigma(i)}+\delta_{\sigma(i)})\Fall{\lam_i+\delta_i}$.

If $a,b$ are nonnegative integers then $a\Fall b=0$ unless
$a\ge b$. So the $\sigma$\_summand vanishes unless
$\mu_{\sigma(i)}+\delta_{\sigma(i)}\ge \lam_i+\delta_i$ for all $i$. 
Summing over $i$, we observe that $|\mu|\le|\lam|$ forces equality for 
each $i$, which implies $\sig(\mu+\delta)=\lam+\delta$. But this
is not possible for $\mu\ne\lam$. \qed


Finally we consider the analogue of the complete symmetric functions,
i.e., $P_d^{r\delta}$ where $d$ stands for $(d,0,\ldots,0)$.

\Proposition P6. For $d\ge0$ we have
$$
P_d^{r\delta}=
{-r\choose d}^{-1}\Sum_{i_j}\prod_{j=1}^n
\left[{-r\choose i_{j-1}-i_j}(x_j-r\delta_j-i_j)\Fall{i_{j-1}-i_j}\right]
$$
where the sum runs through all integer sequences
$d=i_0\ge i_1\ge\ldots\ge i_{n-1}\ge i_n=0$.


\Proof: Let $p_d$ denote the right hand side. Obviously, it has the
right degree $d$ and the coefficient of $x_1^d$ is one. Next we show
that the vanishing condition holds. For this let $x=\mu+r\delta$ with
$|\mu|\le|\lam|$ and $\mu\ne\lam$. Then every summand of $p_d$ is a
multiple of $y_1(y_2-1)\ldots(y_d-d+1)$ where $y_1=\ldots
y_{i_{n-1}}=x_n-r\delta_n=\mu_n$,
$y_{i_{n-1}+1}=\ldots=y_{i_{n-2}}=\mu_{n-1}$ etc. In particular, the
$y_i$ are integers with $0\le y_1\le\ldots\le y_d\le\mu_1$. Now assume
that the product does not vanish, i.e., $y_i\ne i-1$ for all $i$. Then
we claim $y_i\ge i$ for all $i$. Indeed, $y_i\ge y_{i-1}\ge i-1$ and
$y_i\ne i-1$ imply $y_i\ge i$. In particular, $\mu_1\ge y_d\ge d$. But
this is not possible for our choice of $\mu$. This shows $p_d(x)=0$.

Finally, we have to prove symmetry. We are considering the case $n=2$
first. For this we need two basic facts about falling factorials:
(1) $x\Fall{a}(x-a)\Fall{b}=x\Fall{a+b}$ which is obvious and the
Vandermonde identity (2) $(x+y)\Fall{n}=
\sum_{i=0}^n{n\choose i} x\Fall{i}y\Fall{n-i}$.
Letting $i_0=d\ge i_1=i\ge i_2=0$ we obtain that $p_d$ is a multiple of
$$
\sum_i{-r\choose d-i}(x_1-r-i)\Fall{d-i}{-r\choose i}x_2\Fall i.
$$
Applying identity (2) this becomes
$$
\sum_{i,j}{(d-i)!(-r)\Fall{d-i}(-r)\Fall i(-r-i)\Fall{d-i-j}
\over j!(d-i-j)!(d-i)!i!}x_1\Fall j x_2\Fall i.
$$
Using (1), the coefficient becomes $(-r)\Fall{d-i}(-r)\Fall{d-j}\big/
j!(d-i-j)!i!$, which implies symmetry for $p_d(x_1,x_2)$.

Now suppose that $n\ge 3$. Summing over $i=i_{n-1}$ first, we obtain
$$
p_d(x)={-r\choose d}^{-1}\sum_{i=0}^d{-r\choose d-i}{-r\choose
i}x_n\Fall i p_{d-i}(x_1-r-i,\ldots,x_{n-1}-r-i)
$$
By induction we conclude that $p_d$ is symmetric in
$\{x_1,\ldots,x_{n-1}\}$. Summing over $i=i_1$ we obtain
$$
p_d(x)={-r\choose d}^{-1}\sum_{i=0}^d{-r\choose d-i}{-r\choose
i}(x_1-r\delta_1-i)\Fall{d-i} p_i(x_2,\ldots,x_n)
$$
which proves symmetry in $\{x_2,\ldots,x_n\}$.  This concludes the
proof.\qed

\beginsection jack. Difference operators and Jack polynomials

In this section we deduce a different characterization of the
polynomials $P_\lam^{r\delta}$ in terms of difference equations.

Let $\eps_i$ be the $i$-th canonical basis vector in $\CC^n$. The
$i$-th {\it shift operator\/} $T_i$ on functions is defined by $T_i
f(x):= f(x-\eps_i)$, and the {\it $i$-th difference operator\/} is
$\nabla_i:=1-T_i$. These operators commute with each other, and
$T_i,\nabla_i$ also commute with multiplication by $x_j$ for $j\ne i$.

\Definition: Let $t$ be an indeterminate. For $1\le i,j\le n$ put
$$
\Delta_{ij}:=(x_i+t)(x_i+r)^{\delta_j}-x_i^{\delta_j+1}T_i,
\quad
\Delta:=\det (\Delta_{ij}),
\quad
\D:= a_\delta(x)^{-1} \Delta.
$$

\medskip\noindent
Since $\Delta_{ij}$ and $\Delta_{kl}$ commute for $i\ne k$, the
determinant $\Delta$ is well defined. Furthermore, it maps symmetric
polynomials to skew-symmetric ones. Hence $\D$ is a well defined
operator acting on the space of symmetric polynomials. We can develop
$$
\D=D_0t^n+D_1t^{n-1}+\ldots+D_n
$$
into a polynomial where $D_i$ is a difference operator of order $i$ and
$D_0=1$.

\Example. For $r=0$ we obtain
$\D=(t+x_1\nabla_1)\ldots(t+x_n\nabla_n)$, hence
$D_i=e_i(x_1\nabla_1,\ldots,x_n\nabla_n)$.

\noindent
We need the following partial order relation on $\ZZ^n$: we say
$\mu\le\lam$ if $\mu_1+\ldots+\mu_i\le\lam_1+\ldots+\lam_i$ for all
$1\le i\le n$. It has the property that $\lambda$ is a partition if
and only if it is maximal among all its permutations.

\Lemma L1. The operator $\D$ is triangular. More precisely,
$$
\D m_\lam\in \prod_i(\lambda_i+r\delta_i+t)m_\lambda+
\sum_{\mu<\lam} \CC[t] m_\mu.
$$
In particular, $\|deg|\D f\le\|deg|f$ for every symmetric polynomial
$f$.

\Proof: The transition matrix between Schur function $s_\lambda$ and
monomial symmetric functions $m_\mu$ is unitriangular. Hence, it
suffices to prove $\D m_\lam\in
\prod_i(\lambda_i+r\delta_i+t)s_\lambda+
\sum_{\mu<\lam} \CC[t] s_\mu$. Now we multiply by $a_\delta$. By
definition, $a_{\lam+\delta}=a_\delta s_\lam$ is the skew
symmetrization of $x^{\lam+\delta}$. Therefore, it suffices to prove
that $\Delta m_\lam$ is a linear combination of monomials $x^\mu$ with
$\mu\le\lam+\delta$ and that the coefficient of $x^{\lam+\delta}$ has
the indicated form.

For this, observe $\Delta_{ij}= x_i^{\delta_j}
(x_i\nabla_i+r\delta_j+t)+$ lower terms in $x_i$, and that 
$x_i\nabla_i(x_i^m)= mx_i^m+$ lower terms. Thus
$$
\Delta_{ij}x_i^m=(m+r\delta_j+t)x_i^{m+\delta_j}+\hbox{ lower terms
in }x_i.
$$
Expanding the determinant defining $\Delta$, we see
that all monomials occuring in $\Delta m_\lambda$ are
of the form $x^\mu$ with $\mu=\sigma(\lambda)+\tau(\delta)-\eta$
where $\sigma$, $\tau$ are permutations and $\eta\in\NN^n$. All these
$\mu$ are $\le\lam+\delta$. Furthermore, $\mu=\lam+\delta$ implies
$\sigma(\lambda)=\lam$, $\tau=1$, and $\eta=0$. In particular, only
the diagonal term contributes to $x^{\lam+\delta}$. Hence, we obtain
$$
\Delta m_\lam\in \prod_i(\lambda_i+r\delta_i+t)x^{\lambda+\rho}+
\sum_{\mu<\lam+\rho} \CC[t]x^\mu.
$$\qed

For $I\subseteq\{1,\ldots,n\}$, put $\eps_I:=\sum_{i\in I}
\eps_i$, and  $T_If:= (\prod_{i\in I}T_i)f=f(x-\eps_I)$. 
Furthermore, we introduce the functions $\phi_I(x):=
\det c_{ij}^I(x)$ where
$$
c_{ij}^I:=\cases{x_i^{\delta_j+1}&for $i\in I$;\cr
(x_i+r)^{\delta_j}& for $i\not\in I$.\cr}
$$
They behave like ``cut-off functions'':

\Lemma Cutoff. Let $r\ne0$ and $\mu$ be a partition. If $\mu-\eps_I$ is
not a partition then $\phi_I(\mu+r\delta)=0$.

\Proof:
Put $x=\mu+r\delta$ and assume $\mu-\eps_I$ is not a partition. Then
there are two cases:

\item{(1)}$\mu_n=0$ and $n\in I$. Then $x_n=0$ and the $n$-th row of
$c^I(x)$ vanishes. Hence $\phi_I(x)=0$.

\item{(2)}There is $i<n$ such that $i\in I$, $i+1\not\in I$, and
$\mu_i=\mu_{i+1}$. In this case $x_i=x_{i+1}+r$ and $c^I$ has two
proportional rows. Hence, again $\phi_I(x)=0$ and the claim is proved.
\qed

\noindent
Now we prove that each $\Plrd$ is an eigenfunction of $\D$. More
precisely:

\Theorem T4. For each partition $\lam$, we have
$$
\D\Plrd=\prod_i(\lambda_i+r\delta_i+t)\Plrd.
$$
In particular, the action of $\D$ on symmetric polynomials is
diagonalizable with distinct eigenvalues.

\Proof: In view of \cite{L1}, it suffices to show that $\D\Plrd$ satisfies
the vanishing condition. We may exclude the case $r=0$ either by direct
computation or by continuity. Since then $a_\delta(\mu+r\delta)\ne0$ for
all partitions $\mu$, we are left with $\Delta(f)$.

We can expand $\Delta$ as follows:
$\Delta=\sum_{I}d_IT_I$, where $d_I=\|det|d_{ij}^I$ and
$$
d_{ij}^I:=\cases{-x_i^{\delta_j+1}&for $i\in I$;\cr
(x_i+t)(x_i+r)^{\delta_j}& for $i\not\in I$.\cr}
$$
Since $d_I$ is a multiple of $\phi_I$, \cite{Cutoff} holds also for it. Let
$\mu$ be a partition with $|\mu|\le|\lam|$, $\mu\ne\lam$. Then 
$\Delta\Plrd(\mu+r\delta)=\sum_{I}d_I(\mu+r\delta)
\Plrd(\mu-\eps_I+r\delta)$. Since $\Plrd$ satisfies the vanishing condition
it follows from \cite{Cutoff} that
$d_I(\mu+r\delta)\Plrd(\mu-\eps_I+r\delta)=0$ for all $I$. This finishes
the proof of the vanishing condition for $\D\Plrd$ and of the Theorem.\qed

Since the $\Plrd$ form also an eigenbasis for $D_1,\ldots,D_n$ we
obtain:

\Corollary CK1. The difference operators $D_1,\ldots,D_n$ commute
pairwise.

\Corollary C2. Every $\Plrd$ has an expansion of the form 
 $m_\lam+\sum_{\mu<\lam} u_{\lam\mu} m_\mu$.

\Proof: \cite{L1} implies that $\D$ preserves the finite 
dimensional space spanned by $\{ m_\mu\mid \mu\le \lam\}$.  Thus, by
the theorem, it has an eigenvector with the above expansion, which by
the lemma has the same eigenvalue as $\Plrd$. So, they are
equal.\qed

Now we can make the connection to the Jack polynomials. First, we recall
their definition: for an indeterminate $t$ consider the differential
operators
$$
\overline\Delta:=\|det|\big(x_i^{\delta_j}
(t+r\delta_j+x_i{\partial\over\partial x_i}\big);\qquad
\overline\cD(t;r):=a_\delta^{-1}\overline\Delta.
$$
These operators were introduced by Sekiguchi \cite{Seki} and Debiard
\cite{Deb}. Macdonald, \cite{Macd1}, uses them to define the Jack
polynomial $P_\lam^{(1/r)}$: it is the unique eigenvector of
$\overline\cD(t;r)$ which is of the form $m_\lam+\sum_{\mu<\lam}a_\mu
m_\lam$.

\Corollary. The top homogeneous component of $\Plrd$ is $P_\lam^{(1/r)}$.

\Proof: Denote this component by $\Pq$. As observed in the proof of \cite{L1}
$\Delta_{ij}= x_i^{\delta_j}(x_i\nabla_i+r\delta_j+ X)+$ lower terms,
and $x_i\nabla_i=x_i{\partial\over\partial x_i}+$ lower terms. Thus
$\D$ acts on $\Pq$
by $a_\delta^{-1}\det(x_i^{\delta_j}(x_i{\partial\over\partial x_i}
+r\delta_j+ t))=\overline\cD(t;r)$. Consequently $\Pq$ is
an eigenfunction of the Sekiguchi-Debiard operator. The
assertion follows from \cite{C2}.\qed

\beginsection Extra. The extra vanishing theorem

\cite{C2} states that $\Plrd$
contains less monomials than it could according to its definition. In this
section we establish a property of $\Plrd$ which is in a way ``dual''
to that: we are going to prove that $\Plrd$ vanishes at more points than
it should by definition.

Recall that $\lam\subset\mu$ means $\lam_i\le\mu_i$ for all $i$,
i.e., the diagrams are contained in each other. Let $\cP$ be the set of
partitions. A subset $S$ of $\cP$ is called {\it closed\/} if
$\lambda\in S$, $\mu\in\cP$ and $\lambda\subset\mu$ implies $\mu\in
S$. For every closed set $S$ we consider the ideal $\cI_S$ of symmetric
polynomials which vanish at all point $\mu+r\delta$ where $\mu$ a
partition which is {\it not\/} in $S$.

\Theorem TI. Let $S\subseteq\cP$ be closed. Then the ideal $\cI_S$
is stable under the action of $\cD(t;r)$.

\Proof: Again, we may exclude $r=0$ by continuity. Then we have to
show that $\Delta(f)(x)=0$ whenever $f\in\cI_S$ and $x=\mu+r\delta$
with $\mu\in\cP\setminus S$. As in the proof of \cite{T4} it suffices to
consider the products $\phi_I(x)f(x-\epsilon_I)$. Assume this does not
vanish. Then $\mu'=\mu-\epsilon_I\in\cP$ with $f(\mu'+r\delta)\ne0$.
But then $\mu'\in S$, and therefore $\mu\in S$ contradicting the choice
of $\mu$.\qed

\noindent Now we can prove the extra vanishing theorem:

\Theorem C3. Let $\lambda$ and $\mu$ be partitions with
$\lam\not\subset\mu$. Then $\Plrd(\mu+\rho)=0$.

\Proof: Consider the closed subset $S$ of all $\mu$ containing
$\lambda$. We have to show $\Plrd\in\cI_S$. Now for generic $r$, there
exist functions in $\cI_S$ which are  {\it non-zero\/} at
$\lam+r\delta$. (For example, the product of  falling factorials
$\prod_{i,j,k} (x_i-r\delta_j)\Fall{\lam_k}$ is such a function).
The ideal $\cI_S$ is $\cD(t;r)$\_stable. Since $\D$ is diagonalizable,
there must be an eigenfunction of $\D$ in $\cI_S$ with this property.
But this function must be a multiple of some $\Pmrd$. Then
$\Pmrd(\lambda+r\delta)\ne0$ implies $|\mu|\le|\lambda|$. Since
$\Pmrd(\mu+r\delta)\ne0$ we have $\lambda\subset\mu$. Hence
$\mu=\lambda$.\qed

\noindent This can be extended:

\Corollary TII. Let $S\subseteq\cP$ be closed. Then
$\cI_S=\mathop{\oplus}\limits_{\lambda\in S}\CC\Plrd$.

\Proof: Since $\cI_S$ is $\cD$\_stable, there must be $S'\subseteq\cP$
with $\cI_S=\oplus_{\lambda\in S'}\CC\Plrd$. Let $\lambda\in S'$. Since
$\Plrd(\lambda+r\delta)\ne0$, it can not be in $\cP\setminus S$. Hence
$S'\subseteq S$. Conversely, let $\lambda\in S$ and assume there is
$\mu\in\cP\setminus S$ with $\Plrd(\mu+r\delta)\ne0$. Then
$\lambda\subset\mu$ by the extra vanishing theorem. Hence $\mu\in S$
which is impossible. This shows $S\subseteq S'$.\qed

\noindent To round this discussion off, let us mention the following

\Proposition PI. Let $\Lambda$ be the ring of symmetric polynomials (in
$n$ variables). Then every $\cD$\_stable ideal of $\Lambda$ is of the
form $\cI_S$ for some closed subset $S$ of
$\cP$.

\Proof: Clearly, every $\cD$\_stable ideal is of the form
$\oplus_{\lambda\in S}\CC\Plrd$. We have to show that $S$ is closed.
For this we need the following weak form of Pieri's rule proved in the
next section: let $e_1=\sum_ix_i$. Expand $e_1\Plrd=\sum_\mu a_\mu
P_\mu^{r\delta}$. Then $a_\mu\ne0$ whenever
$\mu=\lambda+\epsilon_i\in\cP$. This implies
$\mu=\lambda+\epsilon_i\in S$ whenever $\lambda\in S$ and
$\mu\in\cP$ which is equivalent to $S$ being closed.\qed

\beginsection Pieri. The dehomogeneization operators and the Pieri
formula

Both the $\Plrd$ and the Jack polynomials $P_\lambda^{(1/r)}$ form a
basis of the algebra $\Lambda$ of symmetric polynomials. In particular,
there is a linear isomorphism $\Psi:\Lambda\pfeil\Lambda$ which maps 
$P_\lambda^{(1/r)}$ to $\Plrd$. We are going to show that $\Psi$ can
also be described in terms of difference operators.

For this we define the following variant of $\cD$:
$$
\cE\:=a_\delta^{-1}\det[(x_i+r)^{\delta_j} +t x_i^{\delta_j+1}T_i]=
1+\E_1t+\ldots+\E_nt^n.
$$
Let $\Lambda_d\subseteq\Lambda$ be the subspace spanned by all
$\Plrd$ with $|\lambda|=d$. This is also the space of all polynomials
of degree $\le d$ which vanish in all $\mu+r\delta$ with $|\mu|\le
d-1$.

\Lemma E1. We have $\E_k(\Lambda_d)\subseteq\Lambda_{d+k}$. Moreover,
the effect of $\E_k$ on the top homogeneous components
is multiplication by the elementary symmetric function $e_k$.

\Proof: In the notation of section~\cite{jack}, $\E_k$ has the expansion
$\E_k=a_\delta^{-1}\sum_{|I|=k}\phi_IT_I$. Hence
$\E_kf(x)=a_\delta^{-1}(x)\sum_{|I|=k}\phi_I(x)f(x-\epsilon_I)$. Let
$f\in\Lambda_d$ and $\mu$ be a partition with $|\mu|\le d+k-1$ and
$x=\mu+r\delta$. Then we have $\phi_I(x)f(x-\epsilon_I)=0$. This
means $\E_k f\in\Lambda_{d+k}$.  

For the top homogeneous terms, $T_I= 1$ and 
$\phi_I=\prod_{i\in I}x_i$, hence $\E_k$ acts 
like multiplication by $e_k$.\qed

Now we can prove

\Theorem. a) The difference operators $\E_1,\ldots,\E_n$ commute
pairwise.\Par\noindent
b) Let $\psi:\Lambda\pfeil\CC[\E_1,\ldots,\E_n]$ be the isomorphism
with $\psi(e_k)=\E_k$. Then $\Psi(f)=\psi(f)(1)$ (evaluation at $1$)
for all $f\in\Lambda$.

\Proof: Let $\Lambda_{(d)}$ be the space of symmetric homogeneous
polynomials of degree $d$. Then $\Psi:\Lambda_{(d)}\pf\sim\Lambda_d$
and the inverse is given by taking the top homogeneous component. Thus
\cite{E1} implies that the following diagram commutes
$$
\matrix{\Lambda_{(d)}&\pf\Psi&\Lambda_d\cr
\downarrow\Rechts{e_k}&&\downarrow\Rechts{\E_k}\cr
\Lambda_{(d+k)}&\pf\Psi&\Lambda_{d+k}\cr}
$$
Hence $\Psi(e_kf)=\E_k\Psi(f)$ for all $f\in\Lambda$. This shows a).
Let $f(x)=p(e_1,\ldots,e_k)$. Then
$\Psi(f)=\Psi(p(e_k))=p(\E_k)\Psi(1)=\psi(f)(1)$.\qed

As an application of the theory above we give a new proof of the Pieri
rule for Jack polynomials.

At each lattice point $s=(i,j)$ in the diagram of $\lam$, the 
{\it lower\/} and {\it upper\/} hook-lengths are defined by 
$c_\lam(s) =c_\lam(\alpha;s) := \alpha(\lam_i-j)+(\lam'_j-i+1)$, and 
$c'_\lam(s)=c'_\lam(\alpha;s):= \alpha(\lam_i-j+1)+(\lam'_j-i)$.

Let $\mu\subset\lam$. Then $X(\lam/\mu)$ denotes the set of all boxes
$(i,j)\in\lambda$ such that $\mu_i=\lambda_i$ and $\mu_j'<\lambda_j'$.
Then we define
$$
\psi'_{\lam/\mu}(\alpha):= \prod_{s\in X(\lam/\mu)}
{c_\lam(\alpha;s)/c'_\lam(\alpha;s) \over c_\mu (\alpha;s)/
c'_\mu (\alpha;s)}.
$$
The Pieri formula is the following identity:

\Theorem. For every partition $\lambda$ holds $e_k\Pma=\sum_\lam 
\psi'_{\lam/\mu}(\alpha)\Pla$ where $\lambda$ runs over all partitions
of the form $\mu+\epsilon_I$ for some $I\subset \{1,\ldots,n\}$ 
with $|I|=k$, i.e. $\lam-\mu$ is a 
{\rm vertical $k$-strip}.

\Proof: Applying $\Psi$ to both sides, it suffices to prove
$\E_k\Pmrd= \sum_\lam\psi'_{\lam/\mu}({1/r}) \Plrd$, 
summed over 
$\{ \lam\mid\lam-\mu$ is a vertical k-strip$\}$.  
In any case, $\E_k\Pmrd=\sum_\lam a_{\lambda\mu}\Plrd$ where $\lambda$
is a partition of degree $|\mu|+k$. Evaluating at the point
$x=\lambda+r\delta$ and using the expansion of $\cE_k$ we see
$a_{\lambda\mu}\Plrd(\lam+r\delta)=\E_k\Pmrd(x)=
a_\delta(\lam+r\delta)^{-1}\phi_I(\lam+r\delta)\Pmrd(\mu+r\delta)$.
Hence, it remains to  prove the identity
$$
\psi'_{\lam/\mu}(1/r)=
a_\delta(\lam+r\delta)^{-1}\phi_I(\lam+r\delta)
(c_\lam^{r\delta})^{-1}c_\mu^{r\delta}.
$$
We first calculate $c_\lam^{r\delta}/c_\mu^{r\delta}
={r^{|\lam|-|\mu|}c_\lam'/c_\mu'}$.
Let us put $I':=\{i\not\in I\}$, $J:=\{\lam_i\mid i\in I\}$ 
and $J'=\{\lam_i\mid i\in I'\}$, and for simplicity, let us write 
 $c'_\lam(i,j)$ instead of $c'_\lam(1/r;(i,j))$.
Then it is easy to see that 
for $i\in I$, we have
$c_\lam'(i,j+1)=c_\mu'(i,j)$ unless $j\in J'$. Similarly,
for $i\in I'$ $c_\lam'(i,j)=c_\mu'(i,j)$ unless $j\in J$. 
Taking these cancelations into account we get 
$$
{c_\lam^{r\delta}\over c_\mu^{r\delta}}=
{r^{|\lam|}c_\lam'\over r^{|\mu|}c_\mu'}=
r^k \prod_{i\in I} c_\lam'(i,1)
\prod_{i\in I,j\in J'} {c_\lam'(i,j+1)\over c_\mu'(i,j)}
\prod_{i\in I',j\in J} {c_\lam'(i,j)\over c_\mu'(i,j)}.
$$
On the other hand, $a_\delta^{-1}(\lam+r\delta)\phi_I(\lam+r\delta)$
equals
$$
\prod_{i\in I} (\lam_i+r\delta_i)      \prod_{i\in I,k\in I'\atop i<k}
{(\lam_i+r\delta_i)-(\lam_k+r\delta_k+r) \over 
 (\lam_i+r\delta_i)-(\lam_k+r\delta_k)}\prod_{i\in I,k\in I'\atop k<i}
{(\lam_k+r\delta_k+r)-(\lam_i+r\delta_i) \over 
 (\lam_k+r\delta_k)  -(\lam_i+r\delta_i)}
$$ 
Now the set $\{k\in I'\mid\lam_k=0\}$ equals $\{\lam'_1+1,
\lam'_1+2,\ldots,n\}$, and for $j\in J'$, we have 
$\{k\in I'\mid\lam_k=j\}=\{\lam_{j+1}'+1,\lam_{j+1}'+2,\ldots,\mu'_j\}$. 
Thus the first two products, which can be rewritten as
$\prod_{i\in I} (\lam_i+r(n-i))  
 \prod_{i\in I,k\in I', i<k}
{\lam_i-\lam_k+r(k-i-1) \over \lam_i-\lam_k+r(k-i)},
$ 
become after cancelation,
$$
\prod_{i\in I} (\lam_i+r(\lam_1'-i))  
\prod_{i\in I,j\in J'\atop(i,j)\in\mu}
{\lam_i-j+r(\lam_{j+1}'-i) \over 
\lam_i-j+r(\mu_j'-i) }
=r^k \prod_{i\in I} c_\lam'(i,1)
\prod_{i\in I,j\in J'}
{c_\lam'(i,j+1)\over c_\mu'(i,j)}.
$$ 
Finally, for each $j\in J$, the set $\{i\in I\mid\lam_i=j\}$ equals
$\{\mu_j'+1,\mu_j'+2,\ldots,\lam'_j\}$. Thus after
cancelation the third product 
$\prod_{j\in J,k\in I',k<i}
{\lam_k-\lam_i+r(i-k+1) \over 
\lam_k-\lam_i+r(i-k)}$ becomes
$$
\prod_{j\in J,k\in I'\atop (k,j)\in\mu}
{\lam_k-j+r(\lam_j'-k+1) \over 
\lam_k-j+r(\mu_j'-k+1)}
=\prod_{i\in I',j\in J}
{c_\lam(i,j)\over c_\mu(i,j)}.$$
Since $\psi'_{\lam/\mu}(1/r)= \displaystyle{  \prod_{i\in I',j\in J} 
{ c_\lam(i,j)/c'_\lam(i,j) \over c_\mu (i,j)/ c'_\mu (i,j) } }$,
the result follows.\qed

\beginsection. Scholium.

We close with a conjecture on the ``integral'' form of the Jack 
polynomial. In the homogeneous case, this is the function
$\Jla=c_\lam(\alpha)\Pla$. In the inhomogeneous situation,
consider the function:
$$\Jlrd(x):= (-1)^{|\lam|}c_\lam(1/r) \Plrd(-x).$$  
Various computations suggest the following extension of a
conjecture of Macdonald for $J_\lam^\alpha$.
 
\medskip 
\noindent{\bf Conjecture.} Put $\alpha=1/r$, and write
$\Jlrd=\sum_{\mu\le\lam} \alpha^{|\mu|-|\lam|}a_{\lam\mu}
(\alpha) m_\mu$. Then $a_{\lam\mu}$ is a polynomial
in $\alpha$ with positive integral coefficients.

Recently we have proved Macdonald's original conjecture
as well as the integrality part of the above conjecture.
We shall report on these developments elsewhere.

\beginrefs

\L|Abk:BL|Sig:BL|Au:Biedenharn, L.C.; Louck, J.D.|Tit:A new class of
symmetric polynomials defined in terms of tableaux%
|Zs:Advances in Appl. Math.|Bd:10|S:396--438|J:1989||

\L|Abk:Deb|Sig:De|Au:Debiard, A.|Tit:Polyn{\^o}mes de Tch{\'e}bychev et de
Jacobi dans un espace euclidien de dimension $p$|Zs:C.R. Acad. Sc. Paris%
|Bd:296|S:529--532|J:1983||

\L|Abk:HU|Sig:HU|Au:Howe, R., Umeda, T.|Tit:The Capelli
identity, the double commutant theorem, and multiplicity\_free
actions|Zs:Math. Ann.|Bd:290|S:569--619|J:1991||

\Pr|Abk:Macd1|Sig:M1|Au:Macdonald, I.G.|Artikel:Commuting differential
operators and zonal spherical functions|Titel:Algebraic Groups,
Utrecht 1986|Hgr:A.M. Cohen et al. eds.|Reihe:Lecture Notes Math.|Bd:1271%
|Verlag:Springer-Verlag|Ort:Heidelberg|S:189--200|J:1987||

\Pr|Abk:Macd2|Sig:M2|Au:Macdonald, I.G.|Artikel:Schur functions: theme and
variations|Titel:Actes 28-e S{\'e}minaire Lotharingien|Hgr:-|Reihe:Publ.
I.R.M.A. Strasbourg|Bd:498/S-27|Verlag:-|Ort:-|S:5--39|J:1992||

\B|Abk:Macd3|Sig:M3|Au:Macdonald, I.G.|Tit:Symmetric Functions and Hall
Polynomials, 2nd ed|Reihe:-|Verlag:Oxford Univ. Press|Ort:-|J:1995||
  
\L|Abk:Ok|Sig:Ok|Au:Okounkov, A.|Tit:Quantum immanants and
higher Capelli identities|Zs:Preprint|Bd:-|S:23 pages|J:1995||

\L|Abk:Olsh|Sig:Ol|Au:Olshanski, G.|Tit:Quasi\_symmetric functions and
factorial Schur functions|Zs:Preprint|Bd:-|S:21 pages|J:1995||

\Pr|Abk:Sahi|Sig:Sa|Au:Sahi, S.|Artikel:The spectrum of certain invariant
differential operators associated to Hermitian symmetric spaces%
|Titel:Lie theory and geometry|Hgr:J.-L. Brylinski et al. eds.%
|Reihe:Progress Math.|Bd:123|Verlag:Birkh{\"a}user|Ort:Boston%
|S:569--576|J:1994||

\L|Abk:Seki|Sig:Se|Au:Sekiguchi, J.|Tit:Zonal spherical functions on
some symmetric spaces|Zs:Publ. RIMS, Kyoto Univ.|Bd:12|S:455--459|J:1977||

\L|Abk:Stan|Sig:St|Au:Stanley, R.|Tit:Some combinatorial properties
of Jack symmetric functions|Zs:Adv. Math|Bd:77|S:76--115|J:1989||

\endrefs

\bye